\documentclass{ws-jktr}
\usepackage{epsfig}

\begin{document}

\markboth{Slavik Jablan and Radmila Sazdanovi\' c} {Unlinking
Number and $BJ$Unlinking Gap}


\title{UNLINKING NUMBER AND UNLINKING GAP}

\author{SLAVIK JABLAN, RADMILA SAZDANOVI\'C$^\dag $}

\address{The Mathematical Institute\\
Knez Mihailova 35\\
P.O. Box 367\\
11001 Belgrade\\
Serbia\\
jablans@yahoo.com}

\address{The George Washington University$^\dag $\\
Department of Mathematics\\
Monroe Hall of Government\\
2115 G Street, NW
\\Washington, D.C. 20052 \\
USA\\
radmila@gwu.edu}


\maketitle

\begin{center}
\small
 \parbox{4.1in}{\emph{``It is very easy to define a number of knot invariants so long as one is not concerned with
giving algorithms for their computation. For instance,... One can
change each knot projection into the  projection of a circle by
reversing the overcrossing and undercrossing at, say k double
points of the projection. The minimum number m(k) of these
operations, that is, the minimal number of self-piercings, by
which a knot is transformed into a circle, is a natural measure of
knottedness."\\
K.~Reidemeister} \cite{Re}}
\end{center}

\begin{abstract}

Computing unlinking number is usually very difficult and complex
problem, therefore we define $BJ$-unlinking number and recall
Bernhard-Jablan conjecture stating that the classical
unknotting/unlinking number is equal to the $BJ$-unlinking number.
We compute $BJ$-unlinking number for various families of knots and
links for which the unlinking number is unknown. Furthermore, we
define $BJ$-unlinking gap and construct examples of links with
arbitrarily large $BJ$-unlinking gap. Experimental results for
$BJ$-unlinking gap of rational links up to 16 crossings, and all
alternating links up to 12 crossings are obtained using programs
{\it LinKnot} and {\it K2K}. Moreover, we propose families of
rational links with arbitrarily large $BJ$-unlinking gap and
polyhedral links with constant non-trivial $BJ$-unlinking gap.
Computational results suggest existence of families of
non-alternating links with arbitrarily large $BJ$-unlinking gap.
\end{abstract}



\section{Introduction}

The main topic of this paper is Bernhard-Jablan unlinking number
and $BJ$-unlinking gap. Computing unlinking number is usually very
difficult and complex problem, therefore we define $BJ$-unlinking
number which will be computable due to the algorithmic nature of
its definition.

In order to make precise statements we first need to
introduce basic notation.\\
 The term ``link'' will be used for
both knots and links. Accordingly, under the terms containing the
word ``link" we consider both knot and link
properties/invariants.\\

We use Conway notation \cite{Co,KL}, and the following related
symbols denoting rational and pretzel links. Rational link
diagrams given by Conway symbols $a_1 \, a_2 \, \ldots \, a_n$ are
denoted by $C_{[a_1,a_2,\ldots,a_n]}$. $C_{[a_1,a_2,\ldots,a_n]}$
is a standard diagram of a rational link denoted by
$R_{[a_1,a_2,\ldots,a_n]}$ which corresponds to the continued
fraction $a_n+\cfrac{1}
{a_{n-1}+\cfrac{1}{a_{n-2}+\ldots+\cfrac{1}{a_2+\cfrac{1}{a_1}}}}$
. For a detailed explanation see \cite{KL}. The pretzel link with
$n$ columns of $a_i$ half-twists each is denoted by
$P_{(a_1,a_2,\ldots,a_n)}$ (Fig. 1).

\begin{figure}[th]
\centerline{\psfig{file=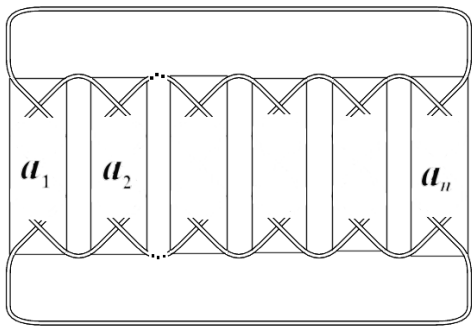,width=2.6in}}
\vspace*{8pt}
\caption{Pretzel link $P_{(a_1,a_2,\ldots,a_n)}$.\label{fig1}}
\end{figure}

The $BJ$-unlinking gap, (Def. 1.2), is motivated by the following
example given by Y.~Nakanishi \cite{Na1} and S.~Bleiler \cite{Bl}.
They have noticed that the rational knot $10_8$, that is
$R_{[5,1,4]}$, has the unlinking number $2$ but its unique minimal
diagram $C_{[5,1,4]}$ (Fig. 2a) has diagram unlinking number $3$.
Unlinking number $2$ can be achieved using non-minimal diagram
illustrated in Fig. 2(b) (crossings needed to be changed are
denoted by circles). The goal of this paper is to determine how
much the unlinking number of any minimal diagram differs from the
unlinking number of a link they represent.

\begin{figure}[th]
\centerline{\psfig{file=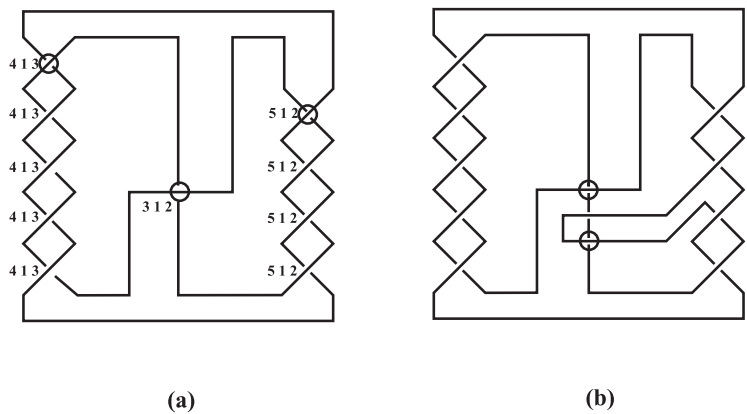,width=4.2in}} \vspace*{8pt}
\caption{The Nakanishi-Bleiler example: (a) the minimal diagram of
the knot $5\,1\,4$ that requires at least three simultaneous
crossing changes to be unknotted. Switching a particular crossing
results in rational knot whose Conway symbol is given next to the
crossing. Notice that all of the diagrams obtained are unknotting
number 2; (b) non-minimal diagram of the knot $5\,1\,4$ from which
we obtain the correct unknotting number
$u(5\,1\,4)=2.$\label{fig2}}
\end{figure}

In this setting we have the following definition:
\begin{definition}
For a given crossing $v$ of a diagram $D$ representing link $L$
let $D_v$ denote the link diagram obtained from $D$ by switching
crossing $v$.

\begin{itemize}
    \item [a)] The {\it unlinking number} $u(D)$ of a link diagram $D$
is the minimal number of crossing changes on the diagram required
to obtain an unlink.
    \item [b)] The \emph{classical unlinking number} of a link $L$,
    denoted by $u(L)$ can be defined by $u(L)=\displaystyle{\min_D} \ u(D)$ where the minimum is
    taken over all diagrams $D$ representing $L$.
        \item [c)] The {\it $BJ$-unlinking number} $u_{BJ}(D)$ of a diagram $D$
is defined recursively in the following manner:
\begin{enumerate}
    \item $u_{BJ}(D)=0$ iff $D$ represents an unlink.
    \item Assume that the sets of diagrams ${\cal D}_k$ with $u_{BJ}(D)\leq k$ are already defined.
     A diagram $D$ has $u_{BJ}(D)= k+1$ if $D \notin {\cal D}_k$ and there exists a crossing $v$
     on the diagram $D$ such that $u_{BJ}(D'_v)=k$ where $D'_v$ is a minimal link diagram
representing the same link as diagram $D_v$ obtained from $D$ by a
crossing change at $v$. Notice that $u_{BJ}(D)$ is well defined
for every diagram $D$ as $D \in {\cal D}_{n!}$ where $n$ is the
number of crossings in $D$.

\end{enumerate}
    \item [d)]$u_M(L)= \displaystyle{\min_D} \ u(D)$ where the minimum is
    taken over all minimal diagrams $D$ representing $L$.
    \item [e)]The {\it $BJ$-unlinking number} $u_{BJ}(L)$ of a link
    $L$ is defined by
$u_{BJ}(L)= \displaystyle{\min_D} \ u_{BJ}(D)$
    where the minimum is taken over all minimal diagrams $D$ representing $L$.
  \end{itemize}

\end{definition}

J.A. Bernhard \cite{Be} in 1994 and independently S.~Jablan
\cite{Ja} in 1995, conjectured that for every link $L$ we have
that $u(L)=u_{BJ}(L)$. In the next section we discuss
$BJ$-conjecture and illustrate the importance of the conjecture on
the example of
pretzel knots whose unlinking number is unknown, except for some small values.\\

\begin{definition}
\begin{itemize}
    \item [a)] The {\it $BJ$-unlinking gap} of a diagram $D$ denoted by $\delta_{BJ}(D)$
is the difference $\delta_{BJ}(D)=u(D)-u_{BJ}(D).$
    \item [b)]  The {\it $BJ$-unlinking gap} of a link $L$ denoted by $\delta_{BJ}(L)$
is defined by $\delta_{BJ}(L)=u_M(L)-u_{BJ}(L)$.
\end{itemize}
\end{definition}

It is natural to consider $\delta'_{BJ}(L)=\displaystyle{\min_D}\
\delta(D)$ where $D$ denotes minimal diagram of $L$. For
alternating links $\delta'_{BJ}(L)=\delta_{BJ}(L)$ and it would be
interesting to check whether this equality holds in greater
generality.

Section 3 contains experimental results-- lists of rational knots
and links up to 16 crossings that have a non-trivial
$BJ$-unlinking gap. Experimental results imply that knots and
links with this property are not so exceptional. In fact, they
represent a considerable portion of rational knots and links,
e.g., about 4$\%$ for $n=15$ or $n=16$.

In the Section 4 we give explicit formulas for $BJ$ unlinking
number, $u_M$ and $BJ$ unlinking gap for several families of knots
(two-bridge and pretzel knots with up to $3$ parameters).
Moreover, we provide formulas for $BJ$-unlinking number and
$BJ$-unlinking gap of the family $R_{[2k,2m,1,2n]}$ and conclude
it has an arbitrarily large unlinking gap. Both Sections 4 and 5
contain experimental results: multi-parameter families of
alternating rational and polyhedral knots and links with positive
$BJ$-unlinking gap. In the Section 5, based on experimental
results, we propose families of non-alternating minimal knot and
link diagrams with possibly arbitrarily large $BJ$-unlinking gap.

\section{Bernhard-Jablan Conjecture}

\bigskip

J.~Bernhard \cite{Be} in 1994 and independently S.~Jablan
\cite{Ja} in 1995, conjectured:
\begin{conjecture}[\textbf{Bernhard-Jablan Conjecture}]
For every link $L$ we have that $u(L)=u_{BJ}(L).$
\end{conjecture}

$BJ$-conjecture holds for all links for which unlinking number is
computed. In particular, it holds for all knots up to 11 crossings
\cite{Liv} and 2-component links up to 9 crossings \cite{Ko2}.
Furthermore, T.~Kanenobu, H.~Murakami and P.~Kohn proved that for
unknotting number 1 rational links unknotting crossing appears in
the minimal diagram \cite{KM,Ko1}.

 Notice that, even if the Bernhard-Jablan Conjecture does not hold for all
links, $BJ$-unlinking number is an upper bound for the unlinking
number.

 We illustrate the importance of the
conjecture by the following example. Consider alternating pretzel
knots $P_{(a,b,c)}$ where $0<a\leq b \leq c$ and $a,b,c$ are all
odd numbers. We show in the Proposition 2.2(b), that
$u_{BJ}(P_{(a,b,c)})=\frac{a+b}{2}$. However, the unknotting
number of these knots is still unknown, except for the smallest
knots, e.g., $P_{(1,3,3)}$ (with unknotting number 2 computed by
W.B.R.~Lickorish \cite{Lic}) and $P_{(3,3,3)}$ (with unknotting
number 3 computed recently by B.~Owens \cite{Ow}).

\begin{figure}[th]
\centerline{\psfig{file=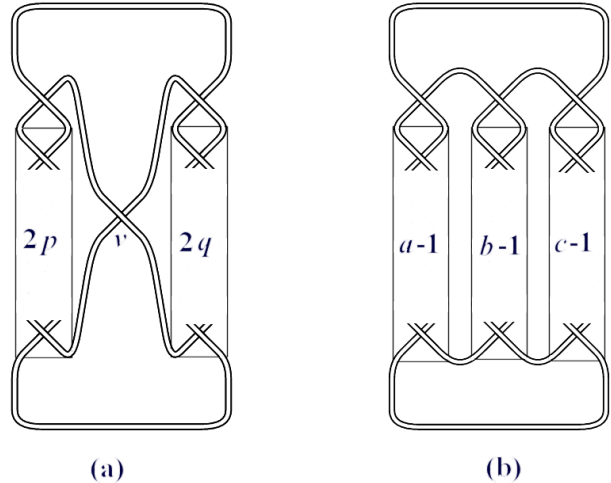,width=3.50in}} \vspace*{8pt}
\caption{(a) Rational knot $R_{[2p-1,1,2q-1]}$; (b) Pretzel knot
$P_{(a,b,c)}$.\label{fig3}}
\end{figure}

\begin{proposition}
\begin{itemize}
    \item [a)] For 2-bridge knots with Conway type $2p-1 \ 1 \ 2q-1$
    denoted by  $R_{[2p-1,1,2q-1]}$
 $($in the fraction form $\frac{4pq-1}{2q}),$
    $p,q>0$ $BJ$-unknotting number satisfies $u_{BJ} (R_{[2p-1,1,2q-1]})=min(p,q).$

    \item [b)]For pretzel knots \footnote{Pretzel knots $P_{(a,b,c)}$ are preserved by
    permutations of symbols a,b,c.} $P_{(a,b,c)}$ with $0<a\leq b \leq c$ and $a,b,c$ are all odd numbers we have that:
    $$u_{BJ}(P_{(a,b,c)})=\frac{a+b}{2}$$
   \end{itemize}
   \end{proposition}
\begin{proof}
\begin{itemize}
    \item [a)]We proceed by induction on $min(p,q)$.
       For $min(p,q)=1$ that is, $p=1$ or $q=1$ our knot is a twist knot
($R_{[1,1,q]}$ or $R_{[p,1,1]}$). Therefore, it has unknotting
number 1, so the proposition holds.
        Assume that proposition holds for $min(p,q) <n$ for
       $n>1$. Since our link is alternating we can work with the
       specific minimal diagram (Fig. 3a). We consider every crossing
       of the diagram. If the crossing is chosen among
       those representing $2p-1$ or $2q-1$ we obtain either
       $R_{[2p-3,1,2q-1]}$ or $R_{[2p-1,1,2q-3]}$, then the inductive step
       is immediate. If we consider the remaining crossing
       $v$, after switching at this crossing and using ambient isotopy
       we get rational knot $R_{[2(p-1),1,2(q-1)]}$. In each step of the inductive
       construction we decrease either $p$ or $q$ or both, therefore we have:
\begin{eqnarray*}
    && u_{BJ} (R_{[2p-1,1,2q-1]})= \\
   &=&1+ min(u_{BJ}
(R_{[2(p-1),1,2q-1]}),u_{BJ} (R_{[2(p-1),1,2(q-1)]}),u_{BJ}
(R_{[2p-1,1,2(q-1)]}))\\
&=&1+min(p-1,q-1) \\
&=& min(p,q)
\end{eqnarray*}
\normalsize
    \item [b)] We proceed by induction on $a+b$.
       For $a+b=2$ our knot $P_{[1,1,c]}$ is a twist knot.
       Such a knot has unknotting number 1, so the proposition holds.
        Assume that proposition holds for $a+b<n$, $n>2$. Since
        our link is alternating we can work with the
       specific minimal diagram (Fig. 3b).
        Assume $a=1$ ($a+b>2$ implies $b,c>1$). We should consider
       every crossing, so we have the following cases:

       \begin{itemize}
        \item If the crossing is chosen among
       those representing  $c$ we obtain
       $P_{[1,b,c-2]}=P_{[1,c-2,b]}$, therefore by inductive
       hypothesis $$u_{BJ}(P_{[1,b,c-2]})=\frac{1+b}{2}$$
        \item If the crossing is chosen among
       those representing $b$ we obtain
       $P_{[1,b-2,c]},$ so by inductive
       hypothesis $u_{BJ}(P_{[1,b-2,c]})=\frac{b-1}{2}$
        \item Crossing change at the remaining crossing
        corresponding to 1 in Conway symbol gives
        the rational knot $R_{[b-2,1,c-2]}$. According to the
      part  a) we have that:
\begin{eqnarray*}
 u_{BJ}(R_{[1,b-2,c-2]})&=&min(\frac{b-1}{2},\frac{c-1}{2})=\frac{b-1}{2}
\end{eqnarray*}
To summarize:
\begin{eqnarray*}
   u_{BJ}(P_{[1,b,c]})&=& 1+ min( u_{BJ}(P_{[1,b,c-2]}), u_{BJ}(P_{[1,b-2,c]}), u_{BJ}(R_{[1,b-2,c-2]}))\\
   &=&1+min(\frac{b-1}{2} ,\frac{1+b}{2})\\
&=&1+ \frac{b-1}{2}=\frac{b+1}{2}
\end{eqnarray*}
which is exactly $\frac{a+b}{2}$ for $a=1$.
       \end{itemize}
\end{itemize}
If $a,b>1$ then we have immediate induction, which completes the
proof. Since, $u_M(R_{[2p-1,1,2q-1]})=u_{BJ}(R_{[2p-1,1,2q-1]})$
and $u_M(P_{[a,b,c]})=u_{BJ}(P_{[1,b,c]})$ the $BJ$-unlinking gap
is equal to zero in both cases.
\end{proof}

Waldhausen \cite{Wa} has proven Smith's conjecture for double
branch covers and we will use it in the following form:

\begin{theorem}[Waldhausen]
Double branched covering $M^{(2)}_K$ of $S^3$ along a knot K is
$S^3$ if and only if K is a trivial knot (unknot).
\end{theorem}

\begin{lemma}[Montesinos]
If $L$ and $m(L)$ are mutant pairs of links then $M^{(2)}_L$ and $
M^{(2)}_{m(L)}$ are homeomorphic $\cite{Mo,Vi}$.
\end{lemma}

\begin{corollary}
Knot $K$ is trivial if and only if $m(K)$ is trivial.
\end{corollary}
\begin{lemma}
If $m(D)$ is a mutation of a diagram $D$ then $u(D)= u(m(D))$.
\end{lemma}
\begin{proof}
Consider an arbitrary 2-tangle $T$ inside a diagram $D$ and a
diagram $m(D)$ obtained from $D$ by mutation of $T$ (the remaining
part of a diagram is intact). If we make $u(D)$ crossing changes
necessary to unknot diagram $D$ and corresponding crossing changes
on a diagram $m(D)$, then, by Corollary 2.5, a diagram obtained
from $m(D)$ also represents an unknot. Therefore, $u(D) \leq
u(m(D))$. Since $D$ can be obtained by mutation on $m(D)$ we have
that $u(D)=u(m(D))$.
\end{proof}

In the following sections we use Tait's first and third
Conjectures. Term minimum diagram stands for minimum crossing
number diagram of a link and reduced diagram is a diagram with no
nugatory crossings.

\begin{theorem}
\begin{itemize}
    \item [a)] $(${\bf Tait's First Conjecture}$)$
A reduced alternating diagram is the minimum diagram of its
alternating link. Moreover, the minimum diagram of a prime
alternating link can only be an alternating diagram. In other
words, a non-alternating diagram can never be the minimum diagram
of a prime alternating link.
    \item [b)] $(${\bf Tait's flyping Conjecture}$)$
   Two reduced alternating diagrams of the same
      link, are related by a finite series of flypes.
\end{itemize}
\end{theorem}

 The first Tait's conjecture \cite{MT1} was proven
independently in 1986 by L.~Kauffman, K.~Murasugi and
M.~Thistlethwaite \cite{Kau,Mu,Th}. The third Tait's conjecture
(Tait's flyping conjecture) was proven by W.~Menasco and
M.~Thistlethwaite \cite{MT1,MT2} and we use it  to prove the
following Corollary:

\begin{corollary}
For every prime alternating link $L$ and its minimal diagrams $D$
and $D'$ the following holds:  $u(D)=u(D')=u_{M}(L)$,
 $u_{BJ}(D)=u_{BJ}(D')=u_{BJ}(L)$, $\delta_{BJ}(D)=\delta_{BJ}(D')=\delta_{BJ}(L).$
\end{corollary}

\begin{proof}
The proof follows from Lemma 2.6 and Tait flyping theorem since
flype can be viewed as a special case od mutation.
\end{proof}

Corollary 2.8. enables us to compute  $BJ$-unlinking number or
$BJ$-unlinking gap using arbitrary minimal diagram of an
alternating link. On the other hand, minimal diagrams of
non-alternating links can have different unlinking numbers, for
example knot $14n_{36750}$
 discovered by A.~Stoimenow \cite{St1}.

\bigskip

\section{Computations of $BJ$-unlinking gap for knots and links}

Experimental results presented in this section are obtained using
{\it Mathematica} based knot theory program {\it LinKnot}
\cite{JS1,JS2,JS3}. For a link given in Conway notation functions
{\bf UnKnotLink} and {\bf fGap} compute $BJ$-unlinking number, the
unlinking number of its fixed minimal diagram, and $BJ$-unlinking
gap $\delta_{BJ} (L)$. Unfortunately, these functions are
dependent on the function {\bf ReductionKnotLink} \cite{Oc} which
sometimes fails in simplifying links. Therefore, for rational
links, we use the {\it LinKnot} function {\bf fGapRat} which is
based on the following Theorem of Schubert \cite{Sch}.
\begin{theorem}[Schubert]
Unoriented rational links $K(\frac{p}{q})$ and $K(\frac{p'}{q'})$
 are ambient isotopic if and only if:
\begin{enumerate}
    \item $p=p'$ and
    \item either $q\equiv q'\pmod p$ or $qq' \equiv 1 \pmod p$
\end{enumerate}
\end{theorem}
The following tables contain Conway symbols of rational knots and
links up to $16$ crossings with a non-trivial $BJ$-unlinking gap,
given according to the number of crossings and whether they are
knots or links. Symbols given in bold denote the links with
$BJ$-unlinking gap 2 (others have $BJ$-unlinking gap 1). The first
column in each table gives the number of crossings, second the
number of knots or links with non-trivial $BJ$-unlinking gap, and
third column their list.

\bigskip

\tiny

\begin{tabular}{|l|l|c|}
  \hline
  $n$ & No. of KL's & List of all KL's \\
  \hline
  $9$ & \begin{tabular}{l} 1 Link \end{tabular}& $4\,1\,4$ \\
  \hline
  $10$ & \begin{tabular}{l}1 Knot \end{tabular} &  $5\,1\,4$ \\
  \hline
  $11$ &\begin{tabular}{l}
    1 Knot \\
    \hline
    4 Links \\

  \end{tabular} &
\begin{tabular}{c}
  $4\,1\,4\,2$ \\
  \hline
  \begin{tabular}{c|c|c|c}
    4 3 4 & 6 1 4 & 5 1 3 2 & 5 1 1 1 3
  \end{tabular} \\
 \end{tabular}

  \\
  \hline
  $12$ & \begin{tabular}{l} 5 Knots \end{tabular} & \begin{tabular}{c|c|c|c|c}
7 1 4 & 5 3 4 & 4 1 4 3 & 6 1 3 2 &  6 1 1 1 3\\
  \end{tabular} \\
  \hline
  $13$ & \begin{tabular}{l}
    7 Knots \\
     \\
    \hline
     \\
    16 Links \\
      \\
      \\

  \end{tabular} &
\begin{tabular}{l}
    \begin{tabular}{c|c|c|c}
      4 4 1 4 & 6 1 4 2 & 4 1 3 1 4 & 5 1 3 2 2 \\
      2 3 1 4 1 2&
    5 1 1 1 3 2   &  5 1 3 1 1 2  &\\

    \end{tabular} \\
    \hline
    \begin{tabular}{c|c|c|c}
{\bf 6 1 6} &   6 3 4 &  8 1 4 &  5 1 5 2
 \\
 5 3 3 2 &6 13 3 & 7 1 3 2 & 3 4 1 3 2
 \\
 4 1 4 2 2 & 5 1 1 1 5 &   6 1 1 2 3 & 7 1 1 1 3
 \\
    2 4 1 3 1 2 & 4 1 1 1 4 2 & 6 1 1 1 2 2 &  4 2 1 1 1 1 3 \\
    \end{tabular}
  \end{tabular}\\
\hline
  $14$&
\begin{tabular}{l}
 31 Knot \\
    \hline
    5 Links \\
 \end{tabular}
   &
   \begin{tabular}{c|c|c|c|c}
7 1 6 & 7 3 4 & 9 1 4 & 4 1 6 3 & 4 3 4 3  \\
5 4 1 4 & 6 1 4 3 &6 1 5 2 & 6 3 3 2 & 7 1 3 3  \\
   8 1 3 2 & 3 3 1 5 2 & 3 5 1 3 2 & 4 1 4 2 3 & 5 1 3 1 4  \\
   5 1 3 2 3 & 5 1 4 2 2 & 6 1 1 2 4 & 6 1 3 2 2 & 7 1 1 2 3  \\
  8 1 1 1 3 & 3 1 3 1 4 2 & 3 1 4 1 3 2 & 3 4 1 3 1 2 & 3 5 1 1 1 3  \\
  5 1 1 1 4 2 & 7 1 1 1 2 2 & 2 1 4 1 3 1 2 & 2 4 1 1 1 3 2 & 4 2 1 1 1 1 4  \\
  5 2 1 1 1 1 3 &  &  & &   \\ \hline
 4 1 4 3 2 & 4 1 6 1 2 & 4 1 3 1 3 2 & 5 1 3 2 1 2 & 5 1 3 1 1 1 2 \\
    \end{tabular}
\\ \hline
\end{tabular}

\bigskip

\normalsize Among the links with $n=13$ crossings we find the
first link $6\,1 \,6$ with the $BJ$-unlinking gap $\delta_{BJ}
=2$.

\smallskip

\footnotesize
\begin{center}
\begin{tabular}{|c|c|c|c|} \hline
  \multicolumn{4}{|c|} {$n=15$} \\ \hline\hline
  \multicolumn{4}{|c|}{ $43$ Knots} \\ \hline
  2 1 3 1 4 1 1 2 &  2 2 4 1 3 1 2 & 2 4 1 1 1 3 1 2 & 2 4 1 1 1 4 2 \\ \hline
   2 4 1 4 2 2 &  3 2 1 1 6 2 & 3 3 1 5 1 2 & 3 3 1 6 2 \\ \hline
    3 4 1 3 2 2  &  4 1 1 1 4 2 2 & 4 1 4 1 1 1 3 & 4 1 4 1 3 2 \\ \hline
     4 1 4 2 2 2  &  4 1 4 3 3  & 4 1 5 1 4 & 4 1 6 1 3 \\ \hline
4 1 7 1 2 & 4 2 1 1 1 1 3 2  & 4 4 3 4 & 4 5 1 1 1 3 \\ \hline
  4 5 1 3 2 &   4 6 1 4 & 5 1 1 1 2 1 4 & 5 1 1 1 3 4 \\ \hline
  5 1 3 1 1 4 & 5 1 3 2 1 3 & 5 1 3 2 4  & 5 1 4 3 2 \\ \hline
   5 1 5 2 2  & 5 1 6 1 2  & 5 3 3 2 2  & 6 1 1 1 2 1 1 2  \\ \hline
    6 1 1 2 3 2 &  6 1 3 1 1 1 2 & 6 1 3 1 4 & 6 1 3 2 1 2  \\ \hline
6 1 3 3 2 & 6 1 4 4 & 6 1 6 2 & 6 4 1 4 \\ \hline
  7 1 1 1 3 2 &  7 1 3 1 1 2  &8 1 4 2   &  \\ \hline

  \multicolumn{4}{|c|}{$63$ Links} \\ \hline
 10 1 4 & 2 1 4 1 3 1 1 2 & 2 3 1 1 1 4 1 2 & 2 3 1 4 1 2 2  \\ \hline
  2 4 1 5 1 2 &  2 5 1 1 1 2 1 2 & 2 5 1 3 2 2 & 3 1 1 1 5 2 2 \\ \hline
   3 4 1 3 1 3 & 3 4 1 5 2  &  3 4 3 3 2 & 3 5 1 1 1 2 2 \\ \hline
    3 5 1 3 3 & 3 6 1 3 2 & 4 1 1 1 1 2 3 2  &  4 1 2 1 1 1 1 1 3 \\ \hline
     4 1 4 1 2 1 2 & 4 1 4 4 2 & 4 1 6 2 2 & 4 2 1 1 3 1 3  \\ \hline
 4 2 1 1 5 2 & 4 2 2 1 1 2 3 & 4 2 3 1 1 1 3 & 4 2 4 1 4 \\ \hline
  4 3 1 1 2 1 3  &  4 3 4 2 2 & 4 4 1 3 1 2 & 5 1 1 1 1 2 4 \\ \hline
   5 1 1 1 2 2 3 & 5 1 1 1 3 2 2  & 5 1 1 1 4 3 & 5 1 3 1 5 \\ \hline
    5 1 5 1 3 & 5 2 1 1 1 1 2 2 & 5 2 1 1 1 2 3  &  5 3 5 2 \\ \hline
     5 4 1 3 2 & 5 5 3 2 & 6 1 1 1 4 2 & 6 1 1 2 5  \\ \hline
 6 1 2 2 4 & 6 1 3 2 3 & 6 1 4 2 2 & 6 1 5 3 \\ \hline
  6 2 1 1 1 1 3  &  6 3 1 2 3 & 6 3 3 3 & {\bf 6 3 6} \\ \hline
   7 1 1 1 2 3 & {\bf 7 1 1 1 5} &  7 1 1 2 2 2 & 7 1 1 3 3 \\ \hline
    7 1 3 1 3 & 7 1 3 4 & {\bf 7 1 5 2}  &  7 3 3 2 \\ \hline
     8 1 1 1 2 2 & 8 1 1 2 3 & 8 1 3 3 & {\bf 8 1 6}  \\ \hline
 8 3 4 & 9 1 1 1 3 & 9 1 3 2 &    \\
   \hline\hline

  \multicolumn{4}{|c|} {$n=16$} \\ \hline\hline
  \multicolumn{4}{|c|}{$138$ Knots} \\ \hline
  10 1 1 1 3  & 10 1 3 2 & 11 1 4 & 2 1 4 1 1 1 3 1 2    \\ \hline
    2 1 4 1 3 1 1 1 2 & 2 1 5 1 1 1 2 1 2  & 2 2 3 1 5 1 2  & 2 3 1 6 2 2 \\ \hline
    2 3 2 1 1 1 1 3 2 &  2 4 1 4 3 2 &  2 4 1 6 1 2 & 2 5 1 1 2 3 2
   \\ \hline
2 6 1 1 2 2 2 & 3 1 1 1 3 1 4 2 & 3 1 2 1 1 1 5 2 & 3 1 2 4 1 3 2
\\ \hline
     3 1 3 1 1 1 4 2 & 3 1 3 1 4 2 2 &  3 1 4 1 1 1 3 2 & 3 1 4 1 3 1 1 2  \\ \hline
     3 1 4 3 3 2 & 3 1 5 1 1 1 2 2 & 3 2 1 4 1 3 2 & 3 2 3 1 5 2\\ \hline
   3 2 4 1 3 1 2 & 3 2 5 1 1 1 3 & 3 3 1 1 1 4 1 2 & 3 3 1 5 1 3\\ \hline
   3 3 3 5 2 & 3 3 4 1 3 2 &  3 4 1 1 1 4 2 & 3 4 1 4 2 2 \\ \hline
   3 4 1 5 1 2 &  3 4 2 1 1 1 1 3 &   3 5 1 1 1 2 1 2 & 3 5 1 3 2 2 \\ \hline
    3 5 1 5 2  & 3 5 3 3 2 & 3 6 1 1 1 2 2  & 3 6 1 1 2 3   \\ \hline
     3 6 1 3 3 & 3 7 1 1 1 3 &  3 7 1 3 2  & 4 1 1 1 4 2 3  \\ \hline
      4 1 2 1 1 1 1 1 4   &  4 1 3 1 4 3 & 4 1 4 1 2 1 3  & 4 1 4 1 4 2 \\ \hline
     4 1 4 2 1 2 2 &  4 1 4 2 2 3  & 4 1 4 2 3 2   & 4 1 4 3 2 2 \\ \hline
      4 1 4 4 1 2 & 4 1 4 4 3 & 4 1 6 1 2 2 &  4 1 6 2 3  \\ \hline
      4 1 8 3  & 4 2 1 1 1 1 3 3 &  4 2 1 1 3 1 4  & 4 2 2 1 1 2 4 \\ \hline
      4 3 1 1 2 1 4 & 4 3 1 6 2  & 4 3 4 2 3 &  4 3 6 3  \\ \hline
      4 4 1 3 1 3  &  5 1 1 1 1 2 3 2  & 5 1 1 1 2 2 4 &  5 1 1 1 4 2 2 \\ \hline
       5 1 2 1 1 1 1 1 3  & 5 1 2 1 1 1 5 & 5 1 4 1 1 1 3 &  5 1 4 1 2 1 2  \\ \hline
      5 1 4 1 3 2 & 5 1 4 4 2  & 5 1 5 1 4  &  5 1 5 2 3  \\ \hline
       5 2 1 1 3 1 3 & 5 2 1 1 5 2 & 5 2 2 1 1 2 3  & 5 2 3 1 1 1 3  \\ \hline
  5 2 3 1 5  & 5 2 4 1 4  & 5 3 1 1 1 5 &  5 3 1 1 2 1 3  \\ \hline
 \end{tabular}

 \begin{tabular}{|c|c|c|c|} \hline
     5 3 3 2 3  &   5 3 4 2 2  & 5 4 1 3 1 2 &  5 4 3 4  \\ \hline
     5 5 1 1 1 3 & 5 5 1 3 2  &   5 6 1 4   &  6 1 1 1 1 2 4 \\ \hline
      6 1 1 1 2 2 3 & 6 1 1 1 3 2 2 & 6 1 1 1 4 3 &  6 1 1 2 3 3   \\ \hline
     6 1 2 2 5 & 6 1 3 1 5  & 6 1 3 2 1 3 &  6 1 3 2 4 \\ \hline
    6 1 3 3 3 & 6 1 4 2 3 & 6 1 4 5  &  6 1 5 1 3 \\ \hline
     6 1 5 2 2 & {\bf 6 1 6 3} &  6 2 1 1 1 1 2 2 &  6 2 1 1 1 2 3  \\ \hline
     6 2 1 1 6  & 6 3 1 2 4  & 6 3 3 2 2  &  6 3 4 3 \\ \hline
         7 1 1 1 4 2 &   7 1 1 3 4 & 7 1 2 2 4 & 7 1 3 1 4  \\ \hline
        7 1 3 2 3   & 7 1 3 3 2 & 7 1 4 2 2 &  7 1 5 3   \\ \hline
   7 2 1 1 1 1 3 &  7 3 1 2 3 &  7 3 3 3  &  7 3 6  \\ \hline
      7 4 1 4  & 8 1 1 1 2 3 & 8 1 1 1 5 & 8 1 1 2 2 2  \\ \hline
      8 1 1 3 3  & 8 1 3 1 3 & 8 1 3 4   & {\bf 8 1 5 2} \\ \hline
     8 3 3 2  &  9 1 1 1 2 2 &   9 1 1 2 3 &   9 1 3 3 \\ \hline
      9 1 6   & 9 3 4  & &   \\ \hline

  \multicolumn{4}{|c|}{$42$ Links} \\ \hline
 2 3 1 4 1 3 2 & 2 3 1 4 3 1 2 & 2 3 1 7 1 2 & 2 3 5 1 3 2 \\ \hline
  2 4 1 3 1 1 1 1 2  & 2 4 1 3 1 2 1 2 & 3 1 1 1 5 3 2 & 3 1 1 1 7 1 2 \\ \hline
   3 4 1 3 1 1 1 2 & 3 4 1 3 2 1 2  &  4 1 3 1 5 2 & 4 1 4 2 2 1 2 \\ \hline
    4 1 4 5 2 & 4 1 5 1 1 1 3 & 4 1 5 2 1 1 2  &  4 1 6 1 4 \\ \hline
     4 1 6 3 2 & 4 1 7 1 1 2 & 4 1 8 1 2 & 4 3 4 1 4  \\ \hline
 4 3 4 3 2 & 4 3 6 1 2 & 5 1 1 1 2 1 3 2 & 5 1 1 1 3 3 2 \\ \hline
  5 1 1 1 5 1 2  &  5 1 3 1 1 1 4 & 5 1 3 1 1 3 2 & 5 1 3 1 2 1 1 2 \\ \hline
   5 1 3 1 3 1 2 & 5 1 3 2 1 4  &  5 1 3 2 3 2 & 5 3 3 1 1 1 2 \\ \hline
    5 3 3 2 1 2 & 6 1 1 1 2 1 1 1 2 & 6 1 1 1 2 2 1 2  &  6 1 3 1 3 2 \\ \hline
     6 1 3 2 1 1 2 & 6 1 3 3 1 2 & 6 1 4 3 2 & 6 1 6 1 2  \\ \hline
 7 1 3 1 1 1 2 & 7 1 3 2 1 2 & &   \\ \hline

 \end{tabular}

\end{center}

\smallskip
\normalsize
\smallskip

First rational knots with the non-trivial unknotting gap
$\delta_{BJ} = 2$ are $6\,1\,6\,3$ and $8\,1\,5\,2$ with $n=16$
crossings. First non-rational alternating knots with
$BJ$-unlinking gap $\delta_{BJ}=1$ appear for $n\ge 12$ crossings:
the pretzel knot $P_{(5,4,3)}$ ($12a_{1242}$) and polyhedral knots
$6^*2.4\,0:3\,0$ ($12a_{970}$), $6^*2.2\,1\,0:4\,0$ ($12a_{76}$),
and $6^*2.2.2.4\,0$ ($12a_{1153}$). In the next section they will
be extended to families (see Def. 4.1) with  $BJ$-unlinking gap
$\delta_{BJ}=1$ .

\section{$BJ$-unlinking gap for Some Families of Alternating Links}    
\medskip

 In this section we explore the effect $2n$-moves \cite{Pr} have on the $BJ$-unlinking number,
$u_M$ and $BJ$-unlinking gap. Applying $2n$-move on an integer
tangle decreases or increases its Conway symbol by $2n$. If we
allow applying $2n$-moves on an arbitrary subset of integer
tangles of a given link we get its infinite families defined
below:

\begin{definition}
For a link or knot $L$ given in an unreduced \footnote{The Conway
notation is called unreduced if in symbols of polyhedral knots or
links elementary tangles 1 in single vertices are included.}
Conway notation $C(L$) denote by $S$ a set of numbers in the
Conway symbol excluding numbers denoting basic polyhedron and
zeros (determining the position of tangles in the vertices of
polyhedron). For $C(L)$ and an arbitrary (non-empty) subset
$\tilde S$ of $S$ the family $F_{\tilde S}(L)$ of knots or links
derived from $L$ is constructed by substituting each $a \in S_f$
by $sgn(a) (|a|+2k_a)$ for $k_a \in N$.
\end{definition}

J.~Bernhard \cite{Be} and D.~Garity \cite{Ga} used this approach
to obtain general formulas for unlinking numbers of the following
families of diagrams of rational knots: $C_{[(2k+1),1,(2k)]}$
($k\ge 2$) and $C_{[(2k+1),(2l+1),(2k)]}$ ($k\ge 2$, $l\ge 0$,
$k>l$) whose  unknotting gap is
$\delta(C_{[(2k+1),(2l+1),(2k)]})=k+l+1-(k+l)=1$. Moreover, the
two-parameter family of rational link diagrams $C_{[2k,1,2l]}$
($k\ge 2$, $l\ge 2$) \cite{Ga} has $u_M(C_{[2k,1,2l]})=k+l-1$ and
$u(C_{[2k,1,2l]})\le l$, so the unlinking gap of a given diagram
is at least $k-1$ and can be made arbitrarily large for a
sufficiently large $k$.

In the similar manner, we try to obtain explicit formulas for
$BJ$-unlinking gap of the infinite family (with up to $k$
parameters) obtained from a link denoted by its Conway symbol
containing $k$ integer tangles. First, we consider rational links
containing only 2 or 3 parameters.

\begin{figure}[th]
\centerline{\psfig{file=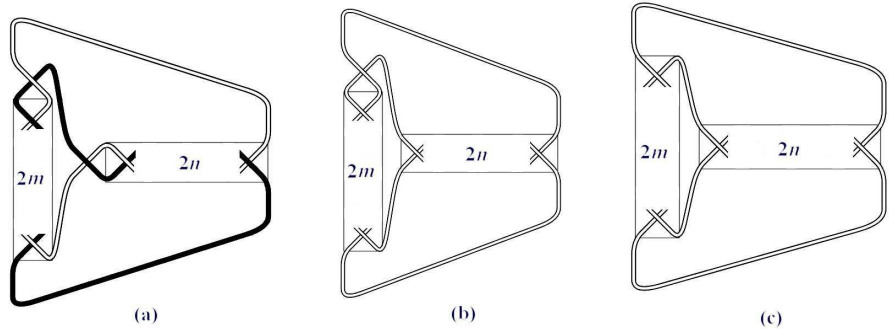,width=4.8in}} \vspace*{8pt}
\caption{(a) Rational link $R_{[2m+1,2n+1]}$; (b) rational knot
$R_{[2m+1,2n]}$; (c) rational knot $R_{[2m,2n]}$.\label{fig6}}
\end{figure}

\begin{lemma}
Let $R_{[a,b]}$ denote a 2-bridge link with the Conway symbol $a \
b$. Then the following holds:
\begin{itemize}
    \item [(a)] If $a,b$ are both odd then for a link $R_{[a,b]}=R_{[2m+1,2n+1]}$
    we have $u_{BJ} (R_{[2m+1,2n+1]})=u_M (R_{[2m+1,2n+1]})=u (R_{[2m+1,2n+1]})=\frac{a+b}{2}=m+n+1$.
    \item [(b)] If $a$ is odd and $b$ is even then for a knot  $R_{[a,b]}=R_{[2m+1,2n]}$ we have
    $u_{BJ} (R_{[2m+1,2n]})=u_M (R_{[2m+1,2n]})=u
    (R_{[2m+1,2n+1]})=n$.
    \item [(c)] If $a,b$ are both even then for a knot
    $R_{[a,b]}=R_{[2m,2n]}$ we have $u_{BJ} (R_{[2m,2n]})=u_M
    (R_{[2m,2n]})=min(m,n)$.
    \footnote{Unknotting number $u(R_{[2m,2n]})$ is an open
    question for the most of knots with $m,n>1.$}
\end{itemize}
\end{lemma}
\begin{proof}
\begin{itemize}
    \item [(a)] Notice that switching $m+1$ crossings corresponding to $2m+1$ and $n$ crossings
    corresponding to $2n+1$, as shown in the picture below, leaves us with
    the unlink of 2 components. Therefore we have: $$u_M (R_{[2m+1,2n+1]})\leq m+n+1$$
    Since every crossing on the diagram is between different
    components with sign equal to 1 belongs to different
    components and has the sign equal to 1, then: $lk
    (R_{[2m+1,2n+1]})=m+n+1$.
Combining these results we have:
\begin{eqnarray*}
   m+n+1&=&lk(R_{[2m+1,2n+1]})\leq  u(R_{[2m+1,2n+1]})\\
        &\leq&  u_{BJ}(R_{[2m+1,2n+1]}) \leq u_M (R_{[2m+1,2n+1]})\\
        &\leq& m+n+1
\end{eqnarray*}
    Therefore $u (R_{[a,b]})=u_{BJ}(R_{[a,b]})=u_M
    (R_{[a,b]})=\frac{a+b}{2}$.

    \item [(b)] In the same manner as in the proof of part a), using signature
    (which is equal to $b$ in this case) instead of linking number we get:
\begin{eqnarray*}
   n&=&\frac{\sigma(R_{[2m+1,2n]}]}{2} \leq  u(R_{[2m+1,2n]})\\
        &\leq&u_{BJ}(R_{[2m+1,2n]})\leq  u_M (R_{[2m+1,2n]})\\
        &\leq& n.
\end{eqnarray*}
 The last inequality is obtained directly from the diagram (switching $n$
crossings corresponding to $2n$ gives the unknot $R_{[2m+1,0]}$)
so we have: $u (R_{[a,b]})=u_{BJ}(R_{[a,b]})=u_M
(R_{[a,b]})=\frac{b}{2}$.

     \item [(c)] We proceed by induction on $min(m,n)$. If $min(m,n)=0$, then we
have an unknot, so the proposition holds. Assume that the
proposition holds for $min(m,n)$ smaller than fixed positive
number $p$. Now consider $min(m,n)=p$. We have a choice of
switching crossing corresponding to $2m$ or $2n$. As a result we
get either $R_{[2(m-1),2n]}$ or $R_{[2m,2(n-1)]}$ with
$min(m-1,n)=min(m,n-1)=p-1$ so according to the induction
hypothesis $u_{BJ}(R_{[2m,2n]})=1+p-1=p=min(m,n)$. Therefore:
$$u(R_{[2m,2n]})\leq u_{BJ}(R_{[2m,2n]})=u_M(R_{[2m,2n]})=min(m,n).$$
    \end{itemize}
\end{proof}

\begin{lemma}
Let $R_{[a,b,c]}$ denote a 2-bridge link with the Conway symbol $a
\ b \ c$. Then the following holds:
\begin{enumerate}

 \item  If $R_{[a,b,c]}=R_{[2k,2l,2m]}$ we have a
2-component link with $u_{BJ} (R_{[2k,2l,2m]})=u_M
(R_{[2k,2l,2m]})=u (R_{[2k,2l,2m]})=k+m$.

\item  If $R_{[a,b,c]}=R_{[2k+1,2l+1,2m+1]}$ then for $(k, m\ge
1)$ has $u_{BJ}(R_{[2k+1,2l+1,2m+1]})=u_M(R_{[2k+1,2l+1,2m+1]})=l+
\min (m,k) +1$.

\item  If $R_{[a,b,c]}=R_{[2k+1,2l,2m+1]}$ $(k, m\ge 1)$ then
$u_{BJ} (R_{[2k,2l+1,2m]})=u_M (R_{[2k,2l+1,2m]})=u
(R_{[2k,2l+1,2m]})=k+l+m+1$.

\item  If $R_{[a,b,c]}=R_{[2k+1,2l,2m]}$ $(k,m\ge 1)$ then $u_{BJ}
(R_{[2k,2l,2m]})=u_M (R_{[2k,2l,2m]})=u (R_{[2k,2l,2m]})=k+m$.

\item  If $R_{[a,b,c]}=R_{[2k,2l+1,2m]}$ $(k,m\ge 1)$ then\\
$u_{BJ}(R_{[2k,2l+1,2m]})=\left\{%
\begin{array}{ll}
    l+max(k,m), & \hbox{$l \leq k,m$;} \\
    k+m, & \hbox{$k,m \leq l$,} \\
\end{array}%
\right.\\$
 i.e., $u_{BJ}(R_{[2k,2l+1,2m]})= \min (k+m, \max(k,m)+l)$
and\\
$u_M(R_{[2k,2l+1,2m]})=\left\{%
\begin{array}{ll}
    k+m, & \hbox{$1 \leq l \leq k,m$;} \\
    k+m-1, & \hbox{$0=l \leq k,m$;} \\
    k+m, & \hbox{$k,m \leq l$.} \\
\end{array}%
\right.$\\

 \item If $R_{[a,b,c]}=R_{[2k+1,2l+1,2m]}$ $(k,m\ge 1)$
then \\
$u_{BJ}(R_{[2k+1,2l+1,2m]})=\left\{%
\begin{array}{ll}
    k+min(l,m), & \hbox{$l,m \leq k$ or $k,m <l$;} \\
    k+l+1, & \hbox{$k,l \leq m$;}
\end{array}%
\right.$
$u_M(R_{[2k+1,2l+1,2m]})=\left\{%
\begin{array}{ll}
    k+min(l,m), & \hbox{$k,m <l$;} \\
    k+min(l+1,m), & \hbox{$l,m \leq k$;} \\
    k+l+1, & \hbox{$k,l \leq m$.}
\end{array}%
\right.$
\end{enumerate}
\end{lemma}

\begin{proof}
Cases $1$, $3$ and $4$ are resolved using of linking number or
signature while the rest require a detailed analysis (similar to
that in the proof of Lemma 4.2) of all possible cases, and will be
omitted.
\end{proof}

\begin{corollary}
\begin{itemize}
    \item [(a)]  Links $R_{[2k,2l+1,2m]}$ have non-trivial gap
if $k,m\ge 2$ and $m\ge l+1$.Then\\
 $\delta_{BJ}(R_{[2k,2l+1,2m]})=\left\{%
\begin{array}{ll}
    min(m,k)-1, & \hbox{$l=0$;} \\
    min(m,k)-l, & \hbox{$l\ge 1$.} \\
\end{array}%
\right.$
    \item [(b)]Family $R_{[2k+1,2l+1,2m]}$ has non-trivial gap
$\delta_{BJ}(R_{[2k+1,2l+1,2m]})= 1$  if $m\ge 2$ and $l+1<m<k$
\footnote{$u_{BJ}(R_{[2k+1,2l+1,2m]})=k+l$}.
\end{itemize}
\end{corollary}

$R_{[a]}$ is a torus knot or link of type $[2,a]$ and therefore
$u_{BJ}=u_M=u$, so both gap and $BJ$-unlinking gap are trivial.
From Lemma 4.2 it follows that all rational links 2 parameters and
$R_{[a,b]}$ have trivial $BJ$-unlinking gap. The same holds for
all 3-parameter families $R_{[a,b,c]}$ except two families listed
in Corollary 4.4.
One can try to extend this classification to rational links with
more parameters and more complicated generating links but the
computations based on parity of parameters and symmetries of the
links are very long and tedious so we give only experimental
results.

All rational links up to 14 crossings with positive $BJ$-unlinking
gap \footnote{Compare with the first table in the Section 3.} are
described by 68 one-parameter\footnote{One parameter family is
obtained by applying the same $2n$-move to all chosen integral
tangles.} families. For all families we predict values of
$BJ$-unlinking number and $BJ$-unlinking gap based on computations
for links with less than 48 crossings. Each family in the
following table is given by its Conway symbol; the next entry is
the number of components followed by experimental results for
$BJ$-unlinking number $u_{BJ}$ and $BJ$-unlinking gap
$\delta_{BJ}.$

\bigskip

\small

\begin{center}

\begin{tabular}{|c|c|c|c|c|}\hline
   & Family & Comp. No. &  $u_{BJ}$ & $\delta_{BJ} $ \\ \hline
  (1) & $(2k+2)\,1\,(2k+2)$ & 2 &  $k+1$  & $k$ \\ \hline
  (2) & $(2k+3)\,1\,(2k+2)$ & 1 &  $k+1$ & 1 \\ \hline
  (3) & $(2k+2)\,1\,(2k+2)\,(2k)$ & 1  & $k+1$  &  $k$ \\ \hline
  (4) & $(2k+2)\,(2k+1)\,(2k+2)$ & 2 & $2k+1$ & $1$ \\ \hline
  (5) & $(2k+4)\,1\,(2k+2)$ & 2  & $k+2$ & $k$ \\ \hline
  (6) & $(2k+3)\,1\,(2k+1)\,(2k)$ & 2 & $k+1$ & $1$ \\ \hline
  (7) & $(2k+3)\,1\,1\,1\,(2k+1)$ & 2  & $k+1$ & $k$ \\ \hline
  (8) & $(2k+5)\,1\,(2k+2)$ & 1  & $k+2$ & $1$ \\ \hline
  (9) & $(2k+3)\,(2k+1)\,(2k+2)$ & 1 &   $2k+1$ & $1$ \\ \hline
  (10) & $(2k+2)\,1\,(2k+2)\,(2k+1)$ & 1 &   $k+1$ & $k$ \\ \hline
  (11) & $(2k+4)\,1\,(2k+1)\,(2k)$ & 1 &   $k+1$ & $1$ \\ \hline
  (12) & $(2k+4)\,1\,1\,1\,(2k+1)$ & 1 &   $k+1$ & $1$ \\ \hline
  (13) & $(2k+2)\,(2k+2)\,1\,(2k+2)$ & 1 &  $k+1$  & $1$ \\ \hline
  (14) & $(2k+4)\,1\,(2k+2)\,(2k)$ & 1 &  $k+2$ & $k$ \\ \hline
  (15) & $(2k+2)\,1\,(2k+1)\,1\,(2k+2)$ & 1  & $k+1$  &  $k$ \\ \hline
  (16) & $(2k+3)\,1\,(2k+1)\,(2k)\,(2k)$ & 1 & $k+1$ & $1$ \\ \hline
  (17) & $(2k)\,(2k+1)\,1\,(2k+2)\,1\,(2k)$ & 1  & $k+1$ & $1$ \\ \hline
  (18) & $(2k+3)\,1\,1\,1\,(2k+1)\,(2k)$ & 1 & $k+1$ & $k$ \\ \hline
  (19) & $(2k+4)\,(2k+1)\,(2k+2)$ & 2  & $2k+2$ & $1$ \\ \hline
  (20) & $(2k+6)\,1\,(2k+2)$ & 2  & $k+3$ & $k$ \\ \hline
  (21) & $(2k+3)\,1\,(2k+3)\,(2k)$ & 2 &   $k+3$ & $1$ \\ \hline
  (22) & $(2k+3)\,(2k+1)\,(2k+1)\,(2k)$ & 2 &   $2k+1$ & $1$ \\ \hline
  (23) & $(2k+4)\,1\,(2k+1)\,(2k+1)$ & 2 &   $k+1$ & $1$ \\ \hline
  (24) & $(2k+5)\,1\,(2k+1)\,(2k)$ & 2 &   $3$ & $1$ if $k=1$ \\
                              & &  &   $k+1$ & $2$ if $k\ge 2$ \\  \hline
  (25) & $(2k+1)\,(2k+2)\,1\,(2k+1)\,(2k)$ & 2  & $2k+1$  &  $1$ \\ \hline
  (26) & $(2k+2)\,1\,(2k+2)\,(2k)\,(2k)$ & 2 & $2k+1$ & $k$ \\ \hline
  (27) & $(2k+3)\,1\,1\,1\,(2k+3)$  & 2  & $k+2$ & $k$ \\ \hline
  (28) & $(2k+4)\,1\,1\,(2k)\,(2k+1)$ & 2 & $k+1$ & $1$ \\ \hline
  (29) & $(2k+5)\,1\,1\,1\,(2k+1)$ & 2  & $k+1$ & $k$ \\ \hline
  (30) & $(2k)\,(2k+2)\,1\,(2k+1)\,1\,(2k)$ & 2  & $k+1$ & $k$ \\ \hline
  (31) & $(2k+2)\,1\,1\,1\,(2k+2)\,(2k)$ & 2 &   $k+1$ & $1$ \\ \hline
  (32) & $(2k+4)\,1\,1\,1\,(2k)\,(2k)$ & 2 &   $k+1$ & $1$ if $k=1,2$\\
       &   &  &   $k$ & $2$ if $k\ge 3$\\ \hline
  (33) & $(2k+2)\,(2k)\,1\,1\,1\,1\,(2k+1)$ & 2 &   $2k$ & $1$ \\ \hline
  (34) & $(2k+5)\,(2k+1)\,(2k+2)$ & 1 &  $2k+2$  & $1$ \\ \hline
  (35) & $(2k+7)\,1\,(2k+2)$ & 1 &  $k+3$ & 1 \\ \hline
  (36) & $(2k+2)\,1\,(2k+4)\,(2k+1)$ & 1  & $k+2$  &  $k$ \\ \hline
  (37) & $(2k+2)\,(2k+1)\,(2k+2)\,(2k+1)$ & 1 & $k+2$ & $1$ \\
  \hline
  (38) & $(2k+3)\,(2k+2)\,1\,(2k+2)$ & 1  & $k+1$ & $k$ \\ \hline
(39) & $(2k+4)\,1\,(2k+2)\,(2k+1)$ & 1 & $k+2$ & $k$ \\ \hline
  (40) & $(2k+4)\,1\,(2k+3)\,(2k)$ & 1  & $k+2$ & $1$ \\ \hline

\end{tabular}

\begin{tabular}{|c|c|c|c|c|}\hline

  (41) & $(2k+4)\,(2k+1)\,(2k+1)\,(2k)$ & 1  & $2k+1$ & $1$ \\ \hline
  (42) & $(2k+5)\,1\,(2k+1)\,(2k+1)$ & 1 &   $k+1$ & $1$ \\ \hline
  (43) & $(2k+6)\,1\,(2k+1)\,(2k)$ & 1 &   $3$ & $1$ if $k=1$ \\
                            & &  &   $k+1$ & $2$ if $k\ge 2$ \\  \hline
  (44) & $(2k+1)\,(2k+1)\,1\,(2k+3)\,(2k)$ & 1 &   $2$ & $1$ for $k=1$\\
                                    & &  &   $k+2$ & $0$ if $k\ge 2$ \\  \hline
  (45) & $(2k+1)\,(2k+3)\,1\,(2k+1)\,(2k)$ & 1 &   $k+1$ & $1$ \\ \hline
  (46) & $(2k+2)\,1\,(2k+2)\,(2k)\,(2k+1)$ & 1 &  $2k+1$  & $k$ \\ \hline
  (47) & $(2k+1)\,1\,(2k+1)\,1\,(2k+2)$ & 1 &  $k+1$ & $k$ \\ \hline
  (48) & $(2k+3)\,1\,(2k+1)\,(2k)\,(2k+1)$ & 1  & $k+1$  &  $1$ \\ \hline
  (49) & $(2k+3)\,1\,(2k+2)\,(2k)\,(2k)$ & 1 & $2k+1$ & $1$ \\ \hline
  (50) & $(2k+4)\,1\,1\,(2k)\,(2k+2)$ & 1  & $k+1$ & $1$ \\ \hline
  (51) & $(2k+4)\,1\,(2k+1)\,(2k)\,(2k)$ & 1 & $2k$ & $1$ \\ \hline
  (52) & $(2k+5)\,1\,1\,(2k)\,(2k+1)$ & 1  & $2$ & $1$ if $k=1$ \\
                              & &  &   $3$ & $0$ if $k\ge 2$ \\  \hline
  (53) & $(2k+6)\,1\,1\,1\,(2k+1)$ & 1  & $k+2$ & $1$ \\ \hline
  (54) & $(2k+1)\,1\,(2k+1)\,1\,(2k+2)\,(2k)$ & 1 &   $k+1$ & $k$ \\ \hline
  (55) & $(2k+1)\,1\,(2k+2)\,1\,(2k+1)\,(2k)$ & 1 &   $k+1$ & $1$ \\ \hline
  (56) & $(2k+1)\,(2k+2)\,1\,(2k+1)\,1\,(2k)$ & 1 &   $k+1$ & $k$ \\ \hline
  (57) & $(2k+1)\,(2k+3)\,1\,1\,1\,(2k+1)$ & 1 &   $$k+1$$ & $k$  \\ \hline

  (58) & $(2k+3)\,1\,1\,1\,(2k+2)\,(2k)$ & 1  & $k+1$  &  $1$ \\ \hline
  (59) & $(2k+5)\,1\,1\,1\,(2k)\,(2k)$ & 1 & $k+1$ & $1$ if $k=1,2$\\
                                 & &  &   $k$ & $2$ if $k\ge 3$ \\  \hline
  (60) & $(2k)\,1\,(2k+2)\,1\,(2k+1)\,1\,(2k)$  & 1  & $2k$ & $1$ if $k=1,2$ \\
                                           & &  &   $2k$ & $k$ if $k\ge 3$ \\  \hline
  (61) & $(2k)\,(2k+2)\,1\,1\,1\,(2k+1)\,(2k)$ & 1 & $2k$ & $1$ \\ \hline
  (62) & $(2k+2)\,(2k)\,1\,1\,1\,1\,1\,(2k+2)$ & 1  & $2k$ & $1$ \\ \hline
  (63) & $(2k+3)\,(2k)\,1\,1\,1\,1\,(2k+1)$ & 1  & $2k$ & $1$ \\ \hline
  (64) & $(2k+2)\,1\,(2k+2)\,(2k+1)\,(2k)$ & 2 &   $2k+1$ & $k$ \\ \hline
  (65) & $(2k+2)\,1\,(2k+4)\,1\,(2k)$ & 2 &   $2k+1$ & $k$ \\ \hline
  (66) & $(2k+2)\,1\,(2k+1)\,1\,(2k+1)\,(2k)$ & 2 &   $2k+1$ & $1$ \\ \hline
  (67) & $(2k+3)\,1\,(2k+1)\,(2k)\,1\,(2k)$ & 2 &   $2k+1$ & $1$\\ \hline
  (68) & $(2k+3)\,1\,(2k+1)\,1\,1\,1\,(2k)$ & 2 &   $2k+1$ & $1$ \\ \hline
\end{tabular}
\medskip

\end{center}

\normalsize

The following results (unless explicitly stated otherwise) are
based on the properties of the generating links and experimental
results for rational, pretzel and polyhedral links up to 16
crossings. First, we present several multi-parameter families of
rational links with an arbitrarily large $BJ$-unlinking gap.

\begin{enumerate}
     \item The family $R_{[(2k),1,(2m),(2n)]}$  has an
     arbitrarily large $BJ$-unlinking gap (see Theorem 4.5).
     \item The family $R_{[(2k+2),1,(2k+2),(2k-1)]}$($k\ge 1$), starting
    with the knot $R_{[6,1,6,3]}$,  has $BJ$-unlinking number $k+1$ and $\delta_{BJ}(R_{[(2k+2),1,(2k+2),(2k-1)]})=k.$
    \item The family $R_{[(4(k+1),1,(2k+3),(2k)]}$($k\ge 1$), starting
    with the knot $R_{[8,1,5,2]}$, has $u_{BJ}(R_{[(4(k+1),1,(2k+3),(2k)]})=k+2$ and $BJ$-unlinking gap is $k+1.$
    \item The family $R_{[(2k+1),1,1,1,(2l+1)]}$, starting with link
    $R_{[5,1,1,1,3]}$, has the unlinking
    gap
    $\delta_{BJ}(R_{[(2k+1),1,1,1,(2l+1)]})=\left\{%
\begin{array}{ll}
    0, & \hbox{$k=l=0$;} \\
    l-1, & \hbox{$k=l>0$;} \\
    l, & \hbox{$k>l>1$.} \\
\end{array}%
\right.$ for $k\geq l\geq0$\\
\item Knots in the family $R_{[(2k),\ldots ,1,\ldots (2k)]}$
($k\ge 2$) and every link in the family $R_{[(2k),\ldots
,(2k),1,(2m),\ldots ,(2m)]}$ ($k,m\ge 2$) have members with
arbitrarily large $BJ$-unlinking gaps. If symbol $k$ occurs $j$
times $BJ$-unlinking number is $[{j\over2}]k$ and
 $\delta_{BJ}=k-1$.
\end{enumerate}

 We use Lemma 4.2 to prove the following theorem about an
example of a family of rational knots with an arbitrarily large
$BJ$-unlinking gap \footnote{Bleiler asked if $\delta(L)=
u_M(L)-u(L)$ has an upper bound \cite{Bl}. Since $\delta_{BJ}(L)
\leq \delta(L)$, Theorem 4.5 provides more examples of links with
unbounded $\delta(L)$ (compare \cite {St2}).}

\begin{theorem}
Let
 $R_{[2k,2m,1,2n]}$ be a rational knot with diagram $C_{[2k,2m,1,2n]}$ $(m,k,n\ge
 0).$
Then the following holds:
\begin{itemize}
    \item [a)] Diagram unlinking number is $u_M(R_{[2k,2m,1,2n]})= n+min(k,m-1)$
    \item [b)]  $BJ$-unlinking number is $$u_{BJ}(R_{[2k,2m,1,2n]})=\left\{
\begin{array}{ll}
    n, & \hbox{if $m \leq n$;} \\
   n+min(k,m-n), & \hbox{if $m > n$.} \\
\end{array}%
\right.$$ \item [c)]  $BJ$-unlinking gap is
$$\delta_{BJ}(R_{[2k,2m,1,2n]})=\left\{
\begin{array}{ll}
    min(k,m-1), & \hbox{if $m \leq n$;} \\
   min(k,m-1)-min(k,m-n), & \hbox{if $m>n$.} \\
\end{array}%
\right.$$
\end{itemize}
\end{theorem}

\begin{proof}
\begin{figure}[th]
\centerline{\psfig{file=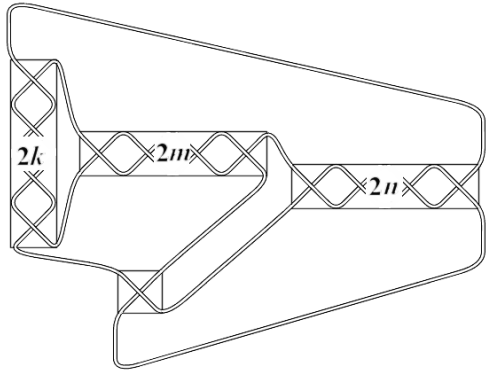, width=2.7in}}
\vspace*{8pt} \caption{Family of rational links $R_{[2k,2m,1,2n]}$
with arbitrarily large $BJ$-unlinking gap. \label{fig5}}
\end{figure}
\begin{itemize}
    \item [a)]
The diagram of $R_{[2k,2m,1,2n]}$ (see Fig. 5) is the minimal as
it is reduced alternating \cite{Kau,Mu,Th}. From Corollary 2.8 it
follows that it is sufficient to consider only one minimal
diagram, so the proof of a) follows from the next lemma:

\begin{lemma}
For $k\geq0, m >0,n>0$
\begin{itemize}
    \item [$($i$)$]$u (C_{[2k,2m,1,2n]})=n+min(k, m-1);$
    \item [$($ii$)$]$u (C_{[2k,2m,-1,2n]})=n-1 +min(k,m-1).$
\end{itemize}
\end{lemma}
\begin{proof} We prove (i) and (ii) simultaneously by induction on
$k+m+n$. For $k+m+n$=2 and $k+m+n=3$ we get: $u (C_{[0,2,1,2]})=u
(C_{[2,1]})=1$, and $u (C_{[0,2,-1,2]})=0$, and $u
(C_{[0,4,1,2]})=u (C_{[2,1]})=1$; $u (C_{[0,4,-1,2]}$ is an
unknot, $u (C_{[0,2,1,4]})=u (C_{[4,1]})=2$, $u (C_{[0,2,-1,4]})=u
(C_{[3]})=1$, $u (C_{[2,2,1,2]})=1$, and $u (C_{[2,2,-1,2]})$ is
an unknot. Assume that the lemma holds for $k+m+n <p$ for $p>3$.
Before we proceed, notice that: $u (R_{[2k,2,1,2]})=u_M
(R_{[2k,2,1,2]})=1$, $u (R_{[2k,2,-1,2]})=u_M
(R_{[2k,2,-1,2]})=0$, and that signature of
$\sigma(R_{[0,2m,1,2n]})=2n$. Therefore, $u (R_{[0,2m,1,2n]})=u_M
(R_{[0,2m,1,2n]})=u(C_{[0,2m,1,2n]})=\frac{\sigma(R_{[0,2m,1,2n]})}{2}=n$
(see Fig. 6). Hence, when $k=0,m=n=2$ lemma holds.

\begin{figure}[th]
\centerline{\psfig{file=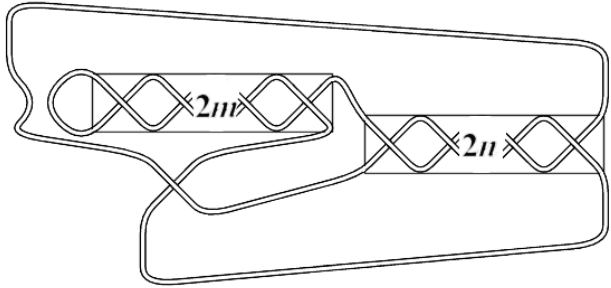, width=3.6in}}
\vspace*{8pt} \caption{ Family of rational links $R_{[0,2m,1,2n]}$
with unlinking number $n$ which can be obtained from a minimal
projection. \label{fig6}}
\end{figure}

\begin{eqnarray*}
u (C_{[2k,2m,1,2n]})&=& 1+ min(\begin{array}{l}
   u (C_{[2(k-1),2m,1,2n]}),u (C_{[2k,2(m-1),1,2n]})\\
   u (C_{[2k,2m,-1,2n]}),u (C_{[2k,2m,1,2(n-1)]}))\\
\end{array}  \\
&=& 1+ min(\begin{array}{l}
  n+min(k-1,m-1),n+min(k,m-1),\\
  u (C_{[2k,2m,-1,2n]}),n-1+min(k,m-1))\\
\end{array}\\
&=& 1+min(n-1+min(k,m-1),u (C_{[2k,2m,-1,2n]}))\\
\end{eqnarray*}

\begin{eqnarray*}
u (C_{[2k,2m,-1,2n]})&=& 1+ min(\begin{array}{l}
   u (C_{[2(k-1),2m,-1,2n]}),u (C_{[2k,2(m-1),-1,2n]})\\
   u (C_{[2k,2m,1,2n]}),u (C_{[2k,2m,-1,2(n-1)]}))\\
\end{array}  \\
&=& 1+ min(\begin{array}{l}
  n-1+min(k-1,m-1),n-1+min(k,m-2),\\
  u (C_{[2k,2m,1,2n]}),n-2+min(k,m-1))\\
\end{array}\\
&=& 1+min(n-2+min(k,m-1),u (C_{[2k,2m,1,2n]}))\\
\end{eqnarray*}
From the equations above we get:
\begin{eqnarray*}
u (C_{[2k,2m,1,2n]})&=& 1+min(n-1+min(k,m-1),u (C_{[2k,2m,-1,2n]}))\\
                    &=&
                    min(n+min(k,m-1),min(n+min(k,m-1),2+u(C_{[2k,2m,1,2n]})))\\
                    &=&min(n+min(k,m-1),2+(C_{[2k,2m,1,2n]}))\\
                    &=&n+ min(k,m-1)\\
                    \end{eqnarray*}
Furthermore:
\begin{eqnarray*} u (C_{[2k,2m,-1,2n]})&= &1+min(n-2+min(k,m-1),u (C_{[2k,2m,1,2n]}))\\
&=& 1+min(n-2+min(k,m-1),min(n-2+min(k,m-1))\\
&=&1+n-2+min(k,m-1)=n-1+min(k,m-1)
\end{eqnarray*}
which completes the proof of the Lemma 4.4 and part a) of the
theorem, i.e., $u_M(R_{[2k,2m,1,2n]})= n+min(k,m-1)$.


\begin{figure}[th]
\centerline{\psfig{file=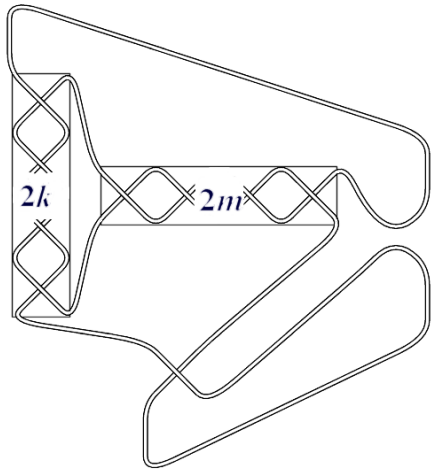, width=2.5In}}
\vspace*{8pt} \caption{Family of rational links
$R_{[2k,2m,1,0]}=R_{[2k,2m]}$ with $BJ$-unlinking number
$\min(m,k).$
 \label{fig5}}
\end{figure}

    \item [b)]
 We proceed by induction on $m+n+k$.
 Using the similar arguments as in Proposition 2.2 and results of the previous lemma
 we get that the proposition holds for:
 \begin{itemize}
    \item $u_{BJ}(R_{[0,2,1,2n]})=u_{BJ}(R_{[1,2n]})=n$
    \item $u_{BJ}(R_{[2k,2m,1,0]})=u_{BJ}(R_{[2k,2m]})=min(k,m)$
    (see Fig.7)
    \item $u_{BJ}(R_{[2k,2,1,0]})=u_{BJ}(R_{[2k,2]})=1$
    \item $u_{BJ}(R_{[2k,0,1,2]})=u_{BJ}(R_{[2k+1,2]})=1$
 \end{itemize}
 In particular, the proposition holds for $m+n+k\leq 3$.
Assume that proposition holds for $m+n+k<p$ for $p \geq 3$. Since
our link is alternating we can work with the specific minimal
diagram (Fig. 3). In the unlinking process we can distinguish 4
cases based on where we perform the crossing change:

  \begin{itemize}
    \item Switch at one of the crossings representing $k$ gives:
    $$u_{BJ}(R_{[2(k-1),2m,1,2n]})=\left\{%
\begin{array}{ll}
    n, & \hbox{$m \leq n$;} \\
    n+min(k-1,m-n), & \hbox{$m>n$.} \\
\end{array}%
\right.$$
    \item Switch at one of the crossings representing $m$ gives:
    $$u_{BJ}(R_{[2k,2(m-1),1,2n]})=\left\{%
\begin{array}{ll}
    n, & \hbox{$m \leq n+1$;} \\
    n+min(k,m-1-n), & \hbox{$m>n+1$.} \\
\end{array}%
\right.$$
    \item Switch at one of the crossings representing $n$ gives:
    $$u_{BJ}(R_{[2k,2m,1,2(n-1)]})=\left\{%
\begin{array}{ll}
    n-1, & \hbox{$m \leq n$;} \\
    n-1+min(k,m-n), & \hbox{$m>n$.} \\
\end{array}%
\right.$$
    \item Switching the crossing representing 1 to $-1$ results in $R_{[2k,2m,-1,2n]}=R_{[2k,2(m-1),1,2(n-1)]}$, therefore:
$$u_{BJ}(R_{[2k,2(m-1),1,2(n-1)]})=\left\{%
\begin{array}{ll}
    n-1, & \hbox{$m \leq n-1$;} \\
    n-1+min(k,m-n+1), & \hbox{$m>n-1$.} \\
\end{array}%
\right.$$
  \end{itemize}
To find $ u_{BJ}(R_{[2k,2m,1,2n]})$ we need to take the minimum
over all 4 cases:
\begin{eqnarray*}
 && u_{BJ}(R_{[2k,2m,1,2n]})= \\
   &=&1+ min(u_{BJ}(R_{[2k,2(m-1),1,2n]}),u_{BJ}(R_{[2k,2(m-1),1,2(n-1)]}),\\
   && u_{BJ}(R_{[2k,2m,1,2(n-1)]}),u_{BJ}(R_{[2k,2(m-1),1,2n]}))\\
&=&1+ \left\{%
\begin{array}{ll}
    min(n,n,n-1,n-1), & \hbox{$m<n$;} \\
    min(n,n,n-1,n-1+min(k,m-n-1)), & \hbox{$m=n$;} \\
    \begin{array}{c}
      min(n+min(k-1,m-n),n+min(k,m-n-1), \\
      n-1+min(k,m-n),n-1+min(k,m-n+1)), \\
    \end{array}, & \hbox{$m \ge n$.} \\
\end{array}%
\right.     \\
&=& 1+ \left\{%
\begin{array}{ll}
   n-1, & \hbox{$m<n$;} \\
   n-1, & \hbox{$m=n$;} \\
   n-1+min(k,m-n), & \hbox{$m > n$.} \\
\end{array}%
\right.     \\
&=&\left\{%
\begin{array}{ll}
    n, & \hbox{$m\leq n$;} \\
    n+min(k,m-n), & \hbox{$m \geq n$.} \\
\end{array}%
\right.
\end{eqnarray*}

\item [c)] Follows from parts a) and b).

\end{proof}
\end{itemize}
\end{proof}

Next we consider the family of pretzel knots $P_{(a,b,c)}$. For
the pretzel knots $P_{(2k+1,2l+1,2m+1)}$ ($k\ge l\ge m\ge 1$) we
proved (Proposition 2.2b) that $u_{BJ}(P_{(2k+1,2l+1,2m+1)})=l+m$
and $\delta_{BJ}(P_{(2k+1,2l+1,2m+1)})=0$. For the families of
pretzel $KL$s with three columns we have the following:
\begin{theorem}
\begin{enumerate}
\item $P_{(2k+1,2l+1,2m+1)}$ has
$u_{BJ}(P_{(2k+1,2l+1,2m+1)})=u_M(P_{(2k+1,2l+1,2m+1)})=l+m$ and
$\delta_{BJ}(P_{(2k+1,2l+1,2m+1)})=0$.

 \item $P_{(2k,2l,2m)}$ has $u_{BJ}(P_{(2k,2l,2m)})=u_M(P_{(2k,2l,2m)})=u(P_{(2k,2l,2m)})=k+l+m$, and therefore
 $\delta_{BJ}(P_{(2k,2l,2m)})=0$ \footnote {Notice that the
linking number guarantees that $u_{BJ}=u_M=u$ and
$\delta=\delta_{BJ}$.};

\item For pretzel knots $P_{(2k+1,2l,2m+1)}$ with $(k\ge m\ge 1)$
we have \footnote {Notice that in first two cases, $l=1$ and $k\ge
l>1$, the signature guarantees that $u_{BJ}=u$ and
$\delta=\delta_{BJ}$.}:\\

\begin{tabular}{|c|c|c|c|}
  \hline
 $P_{(2k+1,2l,2m+1)}$ & $u_{BJ}$ & $u_M$ & $\delta_{BJ}$ \\
  \hline
  $l=1$ & $m+k$ & $m+k$ & $0$ \\
  \hline
  $k\ge l>1$ & $m+k$ & $m+k+1$ & $1$ \\
  \hline
  $l>k\ge 1$ & $m+k+1$ & $m+k+1$ & $0$ \\
  \hline
\end{tabular}

\medskip

 \item $P_{(2k,2l+1,2m)}$ $(k\ge m)$ has
$u_{BJ}(P_{(2k,2l+1,2m)})=u(P_{(2k,2l+1,2m)})=k+l$, and gap
$\delta_{BJ}(P_{(2k,2l+1,2m)}) =m-1$ \footnote {Notice that the
signature guarantees that $u_{BJ}=u$ and $\delta=\delta_{BJ}$.}.
\end{enumerate}

\end{theorem}

Pretzel links $P_{(2k,2l+1,2m)}$ ($k\ge m$) are the example of
links where $BJ$-unlinking number and $BJ$-unlinking gap coincide
with unlinking number and unlinking gap (since half signature
equals $u_{BJ}$), but not with $u_M$, so the gap is non-trivial
and grows as we increase parameter $m$
($u_M(P_{(2k,2l+1,2m)})=k+l+m-1$).

\bigskip

\section{Experimental results and speculations about $BJ$-unlinking gap for polyhedral and non-alternating links}
\medskip
\bigskip
In this section we give experimental results, which (combined with
results from Section 3 for rational links) make computations of
the BJ-unlinking gap complete for alternating links up to $12$
crossings. Furthermore, we propose the family of non-alternating
pretzel link diagrams with an arbitrarily large $BJ$-unlinking
gap.

The first alternating algebraic non-rational link with positive
$BJ$-unlinking gap is the pretzel link $4,4,3$ with $11$-crossings
and $\delta_{BJ} = 1$, and the remaining seven links with
12-crossings are given in the following table:

\medskip

\begin{center}

\small

\begin{tabular}{|c|c|c|c|} \hline
  $2\,2\,1\,1\,2,2,2$ &  $4\,1\,1,3,3$ & $2\,1\,1,3\,1,3\,1$ & $5,4,3$  \\ \hline
    $3\,1,3\,1,2\,1+$ & $(2\,1,2\,1\,1\,1)\,(2,2)$ & $(2\,1,2\,2)\,1\,(2,2)$ & \\ \hline
\end{tabular}

\normalsize

\end{center}

\medskip

Polyhedral knots, defined by Conway \cite{Co}, with $n=12$
crossings and positive $BJ$-unlinking gap are given in the table
below. Second column contains the one-parameter families derived
from these knots, followed by the first step of the unknotting
process which reduces them to families of rational, pretzel, or
polyhedral knots \footnote{ The symbol $\approx$ is used to denote
ambient isotopy between two links; for example, in the first row
symbol $\approx$ means that
$6^*(2k).3\,1.-1.3\,0$ is ambient isotopic to
$(2k-1)\,1\,1\,1\,2\,2$ if $k>1$ and $2\,1\,1\,2\,2$ if $k=1$.}:

\bigskip
\footnotesize

\begin{tabular}{|c|c|c|c|c|}
  \hline
  No. & Knot & Family & Reduction & $\delta_{BJ}$ \\
  \hline
   1 & $6^*2.3\,1:3\,0$ &$6^*(2k).3\,1:3\,0$ & $6^*(2k).3\,1.-1.3\,0$ & $1$ \\
   & $u_{BJ}=2$ & $u_{BJ}=k$&  $k>1$:  $\,\,\approx (2k-1)\,1\,1\,1\,2\,2$ &  \\
    &  &  & $k=1$: $\,\,\approx 2\,1\,1\,2\,2$ &  \\
  \hline
   2 & $6^*2.2\,1\,0:4\,0$ & $6^*(2k).2\,1\,0:4\,0$ & $6^*(2k).2\,1\,0:4\,0:-1$ & $1$ \\
   &  & $u_{BJ}=k+1$ &$\approx 3\,1\,(2k-1)\,2\,2$ &  \\
   \hline
  3 & $6^*2.2\,2\,0:3\,0$ & $6^*(2k).2\,2\,0:3\,0$ & $6^*(2k).2\,2\,0:3\,0:-1$ & $1$ \\
   &  & $u_{BJ}=k+1$  & $\approx 2\,1\,(2k-1),2\,1,2$ &  \\
   \hline
   4 & $6^*2.4\,0:3\,0$ & $6^*2.(2k)\,0:3\,0$ ($k\ge 2$) & $6^*2.(2k)\,0:3\,0:-1$ & $1$ \\
   & $u_{BJ}=2$ & $u_{BJ}=k$ & $\approx (2k)\,1\,1\,2$ &  \\
   \hline
  5 & $6^*2.2.3.3\,0$& $6^*(2k).2.3.3\,0$ & $6^*(2k).2.3.3\,0.-1$ & 1 \\
   &  & $u_{BJ}=k+1$ & $\approx 2\,1\,(2k-1),2\,1,2$ &  \\
\hline
6 & $6^*2.2.2.4\,0$ & $6^*2.2.2.(2k)\,0$ ($k\ge 2$) & $6^*2.2.2.(2k)\,0.-1$  & 1 \\
   &$u_{BJ}=2$  &$u_{BJ}=k$  & $\approx (2k-1)\,1\,1\,1\,1\,2$ &  \\
\hline
 7 & $6^*2.2\,0.3.3\,0$ & $6^*(2k).2\,0.3.3\,0$ & $6^*(2k).2.3.3\,0.-1$ & 1 \\
   &  & $u_{BJ}=k+1$ & $\approx 4\,1\,(2k-1)\,1\,2$  or &  \\
   &  &  & $6^*(2k).2.3.3\,0:-1$ &  \\
   &  &  & $\approx 2\,1\,(2k-1),3,2$ &  \\
\hline
8 & $6^*2.(3,3)$ & $6^*(2k).(3,3)$ & $6^*(2k).(3,3)::-1$ & 1\\
    &  & $u_{BJ}=k+1$ &$k\ge 2$: $\,\,\approx 6^*2.(2k-2):2\,0$ &  \\
    &  &  &$k=1$: $\,\,\approx 2\,1\,1\,1\,2$ &  \\
    \hline
  9 & $6^*2.(3,2).2$ & $6^*2.(2k+1,2).2$ & $6^*2.(2k+1,2).2.-1$  & 1 \\
    &  & $u_{BJ}=k+1$ & $\approx 2\,2\,1\,(2k)$  &  \\
    \hline
    10 & $8^*2:2:.3\,0$ & $8^*(2k):2:.3\,0$ & $8^*2:2.-1..3\,0$ & 1 \\
      & $u_{BJ}=2$     &  & $\approx (2k+1) \, 3 \,2$ (see Lemma 4.3.) &  \\
    \hline

\end{tabular}

\normalsize
\bigskip

Moreover, the following $n=12$-crossing links have $BJ$-unlinking
gap $\delta_{BJ}=1$\footnote{In the first 2 rows we give
2-component links and the third row contains 3-component links.}:

\bigskip

\footnotesize

\begin{tabular}{|c|c|c|c|c|} \hline
 $6^*2.2.2:2\,1\,1\,0$  & $6^*2.2.2:2\,1\,1$ & $6^*(2\,1,2\,2)$ & $6^*2.(2,2)\,1\,1$ & $6^*2.(2,2),2\,0$ \\\hline
 $6^*2.2,(2,2)\,0$  & $6^*2.(2,2).2\,1\,0$  & $6^*(2,2).2\,1:2$ & $6^*(2,2)\,1.2:2\,0$ &  \\\hline
 $6^*2\,1\,1:.(2,2)\,0$  & $6^*2\,1\,1:.(2,2)$ & $6^*2.2\,1.2.2\,0:2\,0$ & $8^*2\,1:.2\,0:2\,0$ & \\\hline
\end{tabular}

\normalsize
\bigskip

The question of finding $BJ$-unlinking gap of non-alternating
links is much more difficult because of the lack of classification
of their minimal diagrams.  For a few classes of non-alternating
links partial results can be obtained using the work of
W.B.R.~Lickorish and M.B.~Thistlethwaite \cite{LT}. Unfortunately,
this is not sufficient to find all minimal diagrams corresponding
to non-alternating link families and compute $BJ$-unlinking gap
for non-alternating links.

The following table contains non-alternating $KL$ diagrams with
$n=11$ and $n=12$ crossings. In all cases $\delta_{BJ} = 1$:

\bigskip \footnotesize

\noindent \begin{tabular}{|c|c|c|c|c|c|} \hline

 $n=11$&  $4\,1\,1,3,-2$ &  $3\,2,3,-3$ & $4,4,-3$ &
$.(3,-2).2$ & $.2.(3,-2)$ \\ \hline

$n=12$ & $5,-3\,1,2\,1$ & $-5,3\,1,2\,1$ & $(-4,2\,1)\,(3,2)$ &
$(-4,-2\,1)\,(3,2)$ & $5,-3\,1,2\,1$ \\ \hline

&  $-5,3\,1,2\,1$ & $(-3,-3)\,(3,2\,1)$
& $(3,3)\,(-3,2\,1)$  & $3:2:-4\,0$ & $-3\,0.2.2\,0.3\,0$ \\
\hline
\end{tabular}
\normalsize

\bigskip

Non-alternating minimal diagrams $4\,1\,1,3,-2$  and $3\,2,3,-3$
(Fig. 8) of the non-alternating knots  $11n_{64}$ and $11n_{122}$
\cite{Ho,Liv} have the unknotting gap $\delta_M=1$. These diagrams
can be extended to two-parameter families of minimal diagrams
$(2k+2)\,1\,1,(2l+1),(-2m)$ and $(2k+1)\,2,(2r+1),-3$ representing
Montesionos knots with the diagram unlinking gap $\delta_M=1$.

\begin{figure}[th]
\centerline{\psfig{file=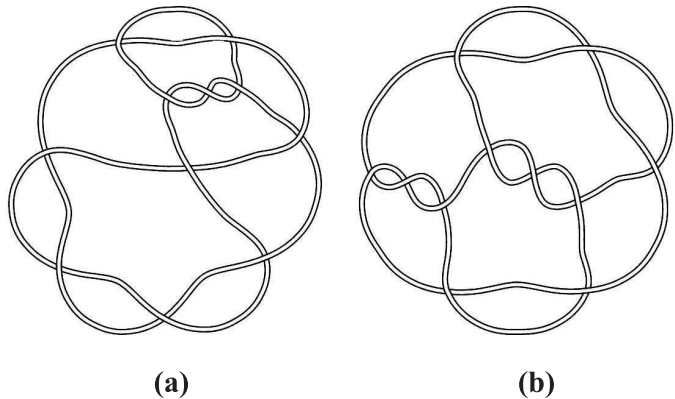,width=3.6in}} \vspace*{8pt}
\caption{The diagrams (a) $4\,1\,1,3,-2$; (b)
$3\,2,3,-3$.\label{fig5}}
\end{figure}

As we described before, even minimal diagrams can have a
non-trivial unlinking gap. Hence, it is not surprising that some
non-minimal diagrams can have a non-trivial unlinking gap.

For example, the 11-crossing non-alternating knot $11n_{138}$
\cite{Ho,Liv} has the non-minimal diagram $3\,1\,1,3,3-$ with the
unknotting gap $\delta_M =1$, while the (fixed) minimal diagram
$3\,1\,1,3,-2\,1$ gives the unknotting number $u=2$ (Fig. 9).

\begin{figure}[th]
\centerline{\psfig{file=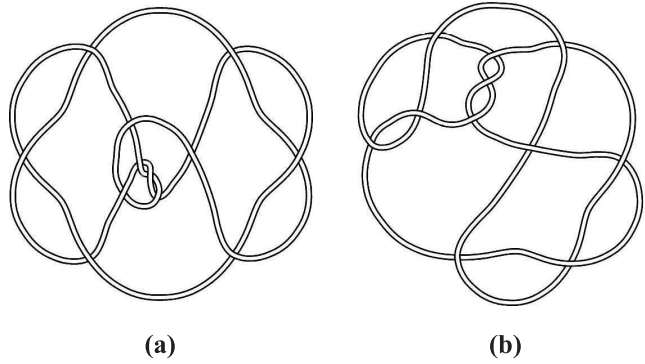,width=3.4in}} \vspace*{8pt}
\caption{(a) Non-minimal diagram $3\,1\,1,3,3-$; (b) the minimal
diagram $3\,1\,1,3,-2\,1$ of the non-alternating knot
$11n_{138}$.\label{fig4}}
\end{figure}

The family of non-alternating pretzel links $P_{(2k,-3,2k)}$
($k\ge 2$) is the candidate for non-alternating link family with
an arbitrarily large unlinking gap (Fig. 10).

\begin{figure}[th]
\centerline{\psfig{file=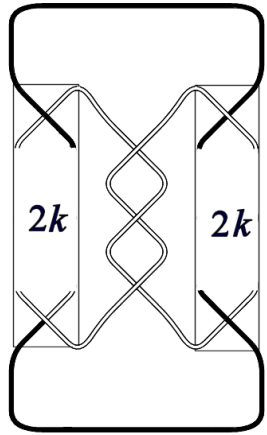,width=2.4in}} \vspace*{8pt}
\caption{The family $P_{(2k,-3,2k)}$ of non alternating minimal
diagrams with an arbitrarily large $BJ$-unlinking
gap.\label{fig10}}
\end{figure}

This family is obtained from the family of rational links
$R_{[2k,1,2k]}=P_{(2k,1,2k)}$ ($k\ge 2$) which is a special case
of the family $R_{[2k,2l+1,2m]}$ from the Corollary 4.4a for
$l=0$, $k=m$ with arbitrarily large $BJ$-unlinking gap
$\delta_{BJ}=k-1$. In the similar manner as in Section 4, we may
obtain that the family of standard diagrams of $P_{(2k,-3,2k)}$
has $BJ$-unlinking number $k$. Furthermore, the unlinking number
of the standard diagram of $P_{(2k,-3,2k)}$ is equal to $2k-1$,
hence the diagram $BJ$-unlinking gap is $k-1$.

Since the classification of all minimal diagrams of the link
family $P_{(2k,-3,2k)}$ is, up to our knowledge, not yet achieved
we are not able to show that the link family $P_{(2k,-3,2k)}$ has
an arbitrarily large unlinking gap.

\section*{Acknowledgments}

We would like to express our gratitude to J\' ozef Przytycki for
his critical reading of the manuscript, corrections, advice and
suggestions.


\bigskip

\begin{thebibliography}{99}


\bibitem[Be]{Be} J.~A. Bernhard, {\it Unknotting numbers and their minimal knot
diagrams}, J. Knot Theory Ramifications, {\bf 3}, 1 (1994) 1--5.

\bibitem[Bl]{Bl} S.~A. Bleiler, {\it A note on unknotting number}, Math. Proc.
Camb. Phil. Soc., {\bf 96} (1984) 469--471.

\bibitem[Co]{Co} J.~Conway, {\it An enumeration of knots and links and some of
their related properties}, in {\it Computational Problems in
Abstract Algebra}, Proc. Conf. Oxford 1967 (Ed. J. Leech),
329--358, Pergamon Press, New York (1970).

\bibitem[Ga]{Ga} D.~Garity,  {\it Unknotting Numbers are not
Realized in Minimal Projections for a Class of Rational Knots},
Proceedings of the "II Italian-Spanish Congress on General
Topology and its Applications" (Italian) (Trieste, 1999). Rend.
Istit. Mat. Univ. Trieste, {\bf 32} (2001), suppl. 2, 59--72
(2002).

\bibitem[Ho]{Ho} J.~Hoste and M.~Thistlethwaite,
{\it Knotscape}, http://www.math.utk.edu/ $^\sim $morwen/

\bibitem[Ja]{Ja} S.~Jablan, {\it Unknotting number and $\infty $-unknotting
number of a knot}, Filomat, {\bf 12}, 1, (1998) 113--120.

\bibitem[JS1]{JS1} S.~Jablan and R.~Sazdanovi\' c, {\it LinKnot},
\newline
http://www.mi.sanu.ac.yu/vismath/linknot/ (2003).

\bibitem[JS2]{JS2}S.~Jablan, R.~Sazdanovic,
LinKnot- Knot Theory by Computer, Series on Knots and Everything,
Volume 21, World Scientific Publishing Co., to appear.

\bibitem[JS3]{JS3} S.~Jablan and R.~Sazdanovi\' c, {\it LinKnot},
http://math.ict.edu.yu/ (2007).


\bibitem [KM] {KM} T.~Kanenobu, H. Murakami {\it Two-bridge knots with Unknotting Number One},
 Proceedings of the American Mathematical Society, {\bf 98}, 3 (1986) 499--502.

\bibitem [Kau]{Kau} L.H. Kauffman, {\it State Models and the
Jones polynomial}, Topology, {\bf 26} (1987) 395--407.

\bibitem [KL]{KL} L.H.~Kauffman and S.~Lambropoulou {\it On the classification of rational
tangles} Advances in Applied Mathematics, {\bf 33}, 2 (2004)
199--237 (see arXiv: math.GT/0212011).


\bibitem [Ko1] {Ko1} P.~Kohn, {\it Two Bridge Links with Unlinking Number One }, Proceedings of the American Mathematical Society,
{\bf 98}, 4 (1991) 1135--1147.

\bibitem [Ko2] {Ko2} P.~Kohn, {\it Unlinking two component links}, Osaka J.
Math., {\bf 30} (1993) 741--752.


\bibitem[Lic] {Lic} W.B.R.~Lickorish,  {\it The unknotting number of a classical knot}, in
{\it Combinatorial methods in topology and algebraic geometry}
(Rocherster, N.Y., 1982), Vol. 44 of {\it Cont. Math.}, 117--121.

\bibitem[LT] {LT} W.B.R.~Lickorish and M.~B.Thislethwaite,  {\it Some links
with non-trivial polynomials and their crossing-numbers}, Comment.
Math. Helvetici, {\bf 63} (1988) 527--539.

\bibitem[Liv] {Liv} C.~Livingston, {\it Knot Tables},
http://www.indiana.edu/$^\sim $knotinfo/ accessed on June 6, 2007.

\bibitem[Mo] {Mo} J.M.~Montesions, {\it Surgery on links and doublebranch covers of $S^3$},
in {\it Knots, groups and $3$-manifolds}, (Ed. L.P. Neuwrith),
Ann. Math. Studies 84, Princeton Univ. Press (1975) 227--259.

\bibitem[MT1] {MT1} W.~W.~Menasco and M.~B.~Thistlethwaite, {\it The Tait flyping
conjecture},  Bull. Amer. Math. Soc., {\bf 25}, 2 (1991) 403--412.

\bibitem[MT2] {MT2} W.~W.~Menasco and M.~B.~Thistlethwaite, {\it
The classification of alternating links}, Annals of Math., {\bf
138} (1993) 113--171.

\bibitem [Mu] {Mu} K.~Murasugi, {\it Jones Polynomials and Classical Conjectures in Knot
Theory II.} Math. Proc. Cambridge Philos. Soc., {\bf 102} (1987)
317--318.



\bibitem[Na1] {Na1} Y.~Nakanishi, {\it Unknotting numbers and knot diagrams with the
minimum crossings}, Math. Sem. Notes Kobe Univ., {\bf 11} (1983)
257--258.

\bibitem[Na2] {Na2} Y.~Nakanishi, {\it Unknotting number and knot
diagram}, Rev. Mat. Univ. Complut., Madrid, {\bf 9} (1996), 2,
359--366.

\bibitem[Oc] {Oc} N.~Imafuji and M.~Ochiai,{\it Computed aided knot theory using Mathematica and MathLink},
J. Knot Theory Ramifications, {\bf 11}, 6 (2002) 945--954.

\bibitem[Ow] {Ow} B.~Owens, {\it Unknotting information from Heegaard Floer homology},
\newline Accepted for publication in Advances in Mathematics; arxiv.org/math.GT/0506485.

\bibitem[Pr] {Pr} J.~Przytycki, {\it $t_k$ moves on links},
Contemporary Math., Vol. 78, Braids - Proceedings of the Santa
Cruz conference on Artin's braid groups (July 1986), 1988, 615-656
\newline (arxiv.org/math.GT/0606633)
\bibitem[Re]{Re} K. Reidemeister, {\it
Knotentheorie.} Ergebn. Math. Grenzgeb., Bd.1; Berlin:
Springer-Verlag (1932) pg. 25. (English translation: Knot theory,
BSC Associates, Moscow, Idaho, USA, 1983).

\bibitem [Ro] {Ro} D. Rolfsen, {\it Knots and Links}, Publish or Perish, 1976 (second
edition, 1990; third edition, AMS Chelsea Publishing, 2003).

\bibitem [Sch] {Sch} H.~Schubert, {\it Knoten mit zwei Br\" ucken}, Math.
Zeit., {\bf 65} (1956) 133--170.


\bibitem[St1] {St1} A.~Stoimenow,  {\it On unknotting numbers and knot
triviadjacency}, Mathematica Scandinavica, {\bf 94}, 2 (2004)
227--248.


\bibitem[St2] {St2} A.~Stoimenow,  {\it On the unknotting number of minimal diagrams}, Mathematics of
Computation, {\bf 72}, 244 (2003) 2043--2057.

\bibitem [Th] {Th} M.~B.Thislethwaite {\it A spanning tree
expansion for the Jones polynomial}, Topology, {\bf 26} (1987),
297--309.

\bibitem [Vi] {Vi} O.~Ya.Viro {\it Nonprojecting isotopies and knots with homeomorphic coverings}, Zap. Nauchn. Semin. LOMI, {\bf 66},
(1976) 133--147, Russian; English transl. in Journal of
Mathematical Sciences, {\bf 12}, 1, (1979), 86--96.

\bibitem [Wa] {Wa} F. Waldhausen {\it \"{U}ber Involutionen der
3-Sph\"{a}re}, Topology, {\bf 8} (1969) 81--91.

\end{thebibliography}
\end{document}